\documentclass[10pt]{article}

\usepackage{amsmath,amssymb,amscd,theorem,pb-diagram,boxedminipage,epic}

{\theorembodyfont{\slshape}\newtheorem{theorem}{Theorem}}
{\theorembodyfont{\slshape}\newtheorem{proposition}[theorem]{Proposition}}
{\theorembodyfont{\slshape}\newtheorem{lemma}[theorem]{Lemma}}
{\theorembodyfont{\slshape}}
{\theorembodyfont{\rmfamily}\newtheorem{note}[theorem]{Note}}
{\theorembodyfont{\slshape}}
\newcommand{\Fq}{\ensuremath{\mathbb{F}_q}}
\newcommand{\Fp}{\ensuremath{\mathbb{F}_p}}
\renewcommand{\O}{\ensuremath{\mathcal{O}}}
\newcommand{\Ot}{\ensuremath{\widetilde{\mathcal{O}}}}
\newcommand{\G}{\ensuremath{\Gamma}}
\newcommand{\g}{\ensuremath{\gamma}}
\newcommand{\og}{\ensuremath{\bar{\gamma}}}
\newcommand{\Q}{\ensuremath{\mathbb{Q}}}
\newcommand{\Qp}{\ensuremath{\mathbb{Q}_p}}
\newcommand{\Qq}{\ensuremath{\mathbb{Q}_q}}
\newcommand{\Zp}{\ensuremath{\mathbb{Z}_p}}
\newcommand{\Zq}{\ensuremath{\mathbb{Z}_q}}

\newcommand{\Z}{\ensuremath{\mathbb{Z}}}

\newcommand{\F}{\ensuremath{\mathbb{F}}}

\newcommand{\ord}{\ensuremath{\text{ord}}}

\newcommand{\ep}{\hfill $\blacksquare$\vspace{\baselineskip}} % End of proof.

\begin{document}

\title{Quasi-quadratic elliptic curve point counting using rigid cohomology}

\author{Hendrik Hubrechts \footnote{Research Assistant of the Research Foundation - Flanders (FWO - Vlaanderen).} \\
\small{Department of mathematics, Katholieke Universiteit Leuven}\\
\small{Celestijnenlaan 200B, 3001 Leuven (Belgium)}\\
\small{\texttt{Hendrik.Hubrechts@wis.kuleuven.be}}}

\maketitle

\begin{abstract}
We present a deterministic algorithm that computes the zeta function of a nonsupersingular elliptic curve $E$ over a finite field with $p^n$ elements in time quasi-quadratic in $n$. An older algorithm having the same time complexity uses the canonical lift of $E$, whereas our algorithm uses rigid cohomology combined with a deformation approach. An implementation in small odd characteristic turns out to give very good results.  \end{abstract}

%\vspace{\baselineskip}
%\noindent \textbf{AMS (MOS) Subject Classification Codes}: 11G20, 11Y99, 12H25, 14F30, 14G50, 14Q05. 
\section{Introduction}\label{sec:intro}
Elliptic curves are a central research object in mathematics, not only centuries and decades ago, but even today with a lot of important unsolved problems concerning such curves. The most notorious example is of course the conjecture of Birch and Swinnerton-Dyer \cite{BirchSwinnertonDyer}, a solution of which is worth a million dollar \cite{ClayMillenium}. In recent times elliptic curves over finite fields have drawn the attention of cryptographers, as Koblitz \cite{KoblitzEC} and Miller \cite{MillerEC} suggested to exploit the group structure on such curves for creating a trapdoor one way function. The motivation for this proposal is that computing discrete logarithms is considered to be very hard for most elliptic curves, while computing the group operation can be done very fast. A very broad exposition can be found in the book \cite{CohenFrey}. Such one way functions can be used in many cryptographic protocols, as for example Diffie-Hellman key exchange \cite{DiffieHellman} or ElGamal encryption \cite{Elgamal}. An important parameter needed for estimating the security level of these applications is the order of the group involved, in this case hence the order of the elliptic curve. We will further on give a brief overview of the large amount of work that has been done on this subject. For now, we content ourselves with noting that counting the number of points on curves over a field of characteristic 2 and of sizes suitable for cryptography can be accomplished in time (far) less than a second.

\subsection{The zeta function and supersingular curves}\label{ssec:zetafunction}
Let $E$ be an elliptic curve defined over the finite field $\Fq$ with $q$ elements, then we can define its zeta function as follows:
\[Z(T) :=\exp\left(\sum_{k=1}^\infty\frac{\#E(\F_{q^k})}{k}T^k\right),\]
where $\#E(\F_{q^k})$ is the number of $\F_{q^k}$-rational points on $E$ (where $E$ is seen as a projective curve). It is well known that $Z(T)$ is actually a rational function, more precisely
\[Z(T) = \frac{qT^2-tT+1}{(1-T)(1-qT)},\qquad t\in\Z,\ |t|\leq 2\sqrt{q}.\]
A proof of this theorem of Hasse and Weil can be found for example in \cite[\S V.2]{Silverman}. The integer $t$ in the zeta function is called the trace of Frobenius, for reasons that will become clear further on in this paper. It is not hard to see that the number $\#E(\F_{q})$ of $\Fq$-rational points on $E$ is precisely $q+1-t$. We can conclude that counting the number of points on $E$ is equivalent to computing its zeta function or its trace $t$.

Curves for which $t\equiv 0\bmod p$ are called \emph{supersingular}, and in \cite[\S V.4]{Silverman} an easy criterion is given for deciding whether a given curve is supersingular. There are only a few possible values for the trace of a supersingular curve, a list with a proof can be found for example in \cite{Waterhouse}. Note that if we are given the zeta function of $E$ over $\Fq$, it is easy to find the zeta function over extension fields of $\Fq$. Indeed, if we denote with $Z_k(T)$ the numerator of the zeta function of $E$ over $\F_{q^k}$, then $Z_k(T)$ equals the following resultant:
\begin{equation}\label{eq:ZetaOverExtension}Z_k(T) = \text{Res}_X(Z_1(X),X^k-T).\end{equation}

\subsection{Point counting algorithms}\label{ssec:pointcounting}
In the following overview we limit our exposition to elliptic curves over finite fields with $p^n$ elements, where $p$ is a small prime number (e.g.\ $p\leq 7$) and $n$ is the relevant parameter. For the complexity estimates --- which are always meant bitwise --- we use the classical Big-Oh notation $\O$, together with the Soft-Oh notation $\Ot$ as defined in \cite[Definition 25.8]{ModernCompAlg}, which ignores logarithmic factors. Using the above remark we also ignore the dependency on $p$ of the algorithms, being irrelevant for very small primes.

A very nice and complete overview of the history of elliptic curve point counting can be found in chapter 17 of the book \cite{CohenFrey} by Cohen, Frey e.a. The first general algorithm is due to Schoof, and improvements by Elkies and Atkin have led to the well known \textsc{sea} algorithm, which runs in time $\Ot(n^4)$ and requires $\O(n^2)$ memory. It is often called `$\ell$-adic', because it works by computing the trace of Frobenius modulo prime numbers $\ell\neq p$. Having done this for enough small primes $\ell$, this allows one to recover the trace.\\

A different approach was considered by Satoh, who found that $p$-adic methods might be much more efficient for small primes $p$ than the technique of Schoof. Satoh's method is based on the \emph{canonical lift} $\mathcal{E}$ of the curve $E$. Let $\Qq$ be the unramified degree $n$ extension of the $p$-adic field $\Qp$, then $\mathcal{E}$ is defined to be the unique (up to isomorphism) lift of $E$ to $\Qq$ which has an endomorphism ring that is isomorphic to the one of $E$, with the isomorphism given by reduction modulo $p$. The idea is then to approximate the $j$-invariant $J$ of this canonical lift modulo an appropriate power of $p$ and afterwards analyzing the action of the $q$th power Frobenius on the lift in order to compute its trace. In later optimizations of the algorithm two main steps arose. First we have to solve an equation $\psi(J,J^\sigma)=0$ over $\Qq$, where $J$ is congruent modulo $p$ to the $j$-invariant of $E$, and $\sigma :\Qq\to\Qq$ is the Frobenius automorphism. A second step consists of computing the norm $\mathcal{N}_{\Qq/\Qp}$ of an element of $\Qq$. Satoh's original algorithm \cite{Satoh} worked in time $\Ot(n^3)$ and required $\O(n^3)$ memory space. After a lot of improvements by Vercauteren \cite{VercauterenQuadMemSatoh}, the \textsc{agm} of Mestre \cite{MestreAGM}, Satoh, Skjernaa and Taguchi (\textsc{sst}) \cite{sst}, and others, a computation time of $\Ot(n^{2.5})$ and space $\O(n^2)$ was achieved. The fastest method however, working for all finite fields of small characteristic, is the patented algorithm of Harley, as described in his e-mail \cite{HarleyECQuadratic}. It requires time $\Ot(n^2)$ and memory $\O(n^2)$, and does not need any precomputations, in contrast to \textsc{sst}. The basic improvements of Harley are fast ways to compute a good representation of $\Qq$, to solve equations of the kind $aX^\sigma+bX+c=0$ over $\Qq$, and to compute the norm $\mathcal{N}_{\Qq/\Qp}$ of an element of $\Qq$. A complete description can be found in section 3.10 of \cite{VercauterenThesis}.

\subsection{An $\Ot(n^2)$, $\O(n^2)$ algorithm using a rigid lift}\label{ssec:QuadraticAlgoRigidLift}
In this paper we describe a new algorithm, that has the same complexity as Harley's result, but is based on a different approach. In \cite{KedlayaCountingPoints} Kedlaya gave an algorithm to compute the zeta function of a hyperelliptic curve of genus $g$ in odd characteristic in time $\Ot(g^4n^3)$ and space $\O(g^3n^3)$. It uses not the canonical lift (for genus one curves), but a rigid lift, which is trivial to compute. If we take the de Rham cohomology of this lifted curve, a Lefschetz fixed point theorem of Monsky and Washnitzer tells us that the characteristic polynomial of the Frobenius operator on this cohomology yields the zeta function of the curve. Three points are crucial. First, if the lift is well-chosen, we can effectively compute in this Monsky-Washnitzer cohomology due to it being isomorphic to the de Rham cohomology of the algebraic lift. Second, by cutting out Weierstrass points, the action of Frobenius is readily computable. And third, factoring the $q$th power Frobenius in repeated applications of the $p$th power Frobenius makes sure that the appearing power series converge good enough. Later on Denef and Vercauteren extended Kedlaya's method to the technically more difficult case of characteristic 2 in \cite{DenefVercauteren}.

In \cite{LauderDeformation}, Lauder used deformation in order to compute the zeta function of higher dimensional varieties. This works by putting the variety in a well-chosen one parameter family, say with formal parameter $\G$, and computing the general matrix $F(\G)$ of the $p$th power Frobenius. As shown by Dwork in \cite{DworkDeformation} such a matrix satisfies a differential equation, the Picard-Fuchs equation of the deformation, and this equation allows fast recovery of $F(\G)$ modulo a certain power of $\G$. In a next step the matrix $F(\G)$ is specialized to $F(\g)$ for some $\g\in\Qq$, and computing the matrix of the $q$th power Frobenius yields then the zeta function. In \cite{HubrechtsHECOdd} and \cite{HubrechtsHECEven} we followed a suggestion of Lauder to try to combine such a deformation with Kedlaya's and Denef and Vercauteren's algorithm, and this resulted in an $\Ot(n^{2.667})$ algorithm for hyperelliptic curves in certain families. The most time consuming step in these algorithms is the computation of the `norm' of the matrix $F(\g)$, i.e.\ the ordered product of its conjugates. For elliptic curves we show in this paper that all curves can be put in a good family, and that we can reduce the problem to computing the norm of just \emph{one} element of $\Qq$. Using Harley's fast norm computation algorithm this gives then the aforementioned complexities. We note that Harley's other basic improvements are also used in our algorithm.\\

We briefly sketch the structure of this paper. In section \ref{sec:parameterFamily} we describe how to put a general curve in a good linear family defined over the prime field, and in the next two sections we repeat briefly how the theory of \cite{HubrechtsHECOdd} and \cite{HubrechtsHECEven} allows us to compute the matrix of the $p$th power Frobenius for curves in such a family. In addition we explain how to recover an integral matrix of Frobenius, which is not guaranteed by the original algorithms of \cite{HubrechtsHECOdd} and \cite{HubrechtsHECEven}. In the fifth section is shown how to compute the trace of Frobenius from this matrix, and in the last section we present an overview of the algorithm and some results obtained with an implementation of (a variant of) the algorithm.\\

The author wants to thank Jan Denef for his help on the problem of finding an integral matrix of Frobenius in characteristic 2 and Denef and Wouter Castryck for their comments on an early version of this paper. \section{The curve placed in a one parameter family}\label{sec:parameterFamily}
Let $E$ be a nonsupersingular elliptic curve over a finite field $\Fq$, given by its Weierstrass equation. We will show in this section how to reduce efficiently the equation of $E$ to another equation over $\Fq$, defining $E'$, such that this last one can be tackled directly using the deformation technique of sections \ref{sec:FrobInOddChar} and \ref{sec:FrobInChar2}. The resulting elliptic curve $E'$ will be isomorphic to the original curve or to its quadratic twist, which we denote by $\text{Twist}(E)$. It is well known that the trace of Frobenius $t$ of $E$ equals minus the trace of Frobenius of $\text{Twist}(E)$, and hence it suffices to work with $E'$. Note that it will be clear in each case which of the two isomorphisms $E'\cong E$ or $E'\cong \text{Twist}(E)$ holds. We have to stress that these results are certainly not new, but we did not find a good reference, and the explicit way to find the curve $E'$ is an important part of a concrete implementation of the algorithm.

\subsection{Odd characteristic}\label{ssec:OddCharParam}
Let $p$ be an odd prime and $\Fq$ a finite field of order $q=p^n$. We suppose that the elliptic curve $E$ over $\Fq$ is given by
\begin{equation}\label{eq:GeneralOdd}
Y^2=X^3+aX^2+bX+c,\quad a,b,c\in\Fq.\end{equation}
If $p\neq 3$ the translation $X\mapsto X-a/3$ removes the term with $X^2$ in (\ref{eq:GeneralOdd}), so we can suppose in this case that $a=0$. If $c=0$ this can be written as $Y^2=X^3+\og X$ with $\og:= b$, a form suitable for section \ref{sec:FrobInOddChar}, so we may assume that $c\neq 0$. Similarly we can assume that $b\neq 0$. The notation $(\Fq)^2$ will be used for the set of squares of $\Fq$.
\begin{proposition}\label{prop:CharOdd}
Let $\og:=b^3/c^2$ and let $E'$ be the elliptic curve over $\Fq$ defined by $Y^2=X^3+\og X+\og$.
If $b/c\in(\Fq)^2$ we have that $E'\cong E$ (over $\Fq$), and otherwise $E'\cong \text{Twist}(E)$.\end{proposition}
\textsc{Proof.} Let $d$ be a nonsquare in $\Fq$ if $b/c\not\in(\Fq)^2$, and $d:=1$ otherwise. Then there exists $\lambda\in\Fq$ such that $\lambda^2=\frac{b}{cd}$, and the change of variables $Y\mapsto \lambda^{-3}Y$, $X\mapsto \lambda^{-2}X$ transforms $Y^2=X^3+bd^2X+cd^3$ into $Y^2=X^3+({b^3}/{c^2})X+{b^3}/{c^2}$.
Is $d$ is a nonsquare the equation $Y^2=X^3+bd^2X+cd^3$ gives precisely the quadratic twist of $E$.\ep

\noindent Now we take $p=3$.
If $a=0$ in (\ref{eq:GeneralOdd})\footnote{All such curves are in fact supersingular because their $j$-invariant is zero.}, we can again use proposition \ref{prop:CharOdd}, and if $a\neq 0$ the translation $X\mapsto X-\frac{b}{2a}$ removes the term with $X$ in (\ref{eq:GeneralOdd}). So we can suppose for the next proposition that $b=0$ and $a\neq 0$.

\begin{proposition}\label{prop:Char3}
Let $p=3$ and $E$ be given by $Y^2=X^3+aX^2+c$. Define $\og:=c/a^3$ and the elliptic curve $E'$ with equation $Y^2=X^3+X^2+\og$. If $a\in(\Fq)^2$ we have that $E'\cong E$, and otherwise $E'\cong \text{Twist}(E)$.\end{proposition}
\textsc{Proof.} If we `twist' $E$ using $a^{-1}$, we find $Y^2=X^3+X^2+c/a^3$, and now we can finish as in the proof of proposition \ref{prop:CharOdd}.\ep

We can conclude that given any elliptic curve in odd characteristic, we can always find $\og\in\Fq$ and some polynomial $Q(X,\G)$ over $\Fp$ such that the following holds: $Q(X,\G)$ is monic of degree 3 in $X$ and linear in $\G$, and it suffices to compute the zeta function of $Y^2=Q(X,\og)$. In addition, this can be done very fast. Indeed, the complexity is dominated by verifying whether $b/c$ (or $a$) is a square in $\Fq$, and as $x\in(\Fq)^2$ is equivalent to $x^{(q-1)/2}=1$, this can certainly be done in time $\Ot(n^2)$ and space $\O(n)$.

In section \ref{sec:FrobInOddChar} we will need that $Y^2=Q(X,0)$ defines an elliptic curve over $\Fp$, but this can always be achieved by the translation $\G\mapsto\G+\alpha$ for some $\alpha\in\Fp$. It is interesting to make the degree in $\G$ of the resultant $\text{Res}_X(Q(X,\G);\frac{\partial}{\partial X}Q(X,\G))$ as small as possible (where we interpret $Q(X,\G)\in\Zp[X,\G]$ for the moment). In proposition \ref{prop:CharOdd} this will be 3 and in proposition \ref{prop:Char3} we find degree 2. If $\og\in (\Fq)^2$ in proposition \ref{prop:CharOdd}, we can twist over $1/\sqrt{\og}$ and find $Y^2=X^3+X+\og'$ for some $\og'\in\Fq$, which also gives a second degree resultant. Although this requires the computation of a square root in $\Fq$, it might still be advantageous in the end.

\subsection{Characteristic 2}\label{ssec:EvenCharParam}
We now take $q=2^n$ and $E$ a nonsupersingular curve over $\Fq$ given by
\[Y^2+a(X+b)Y=X^3+cX^2+dX+e\quad\text{with}\quad a,b,c,d,e\in\Fq.\] The fact that $E$ is not supersingular is easily seen to be equivalent to $a\neq 0$.
The translation $X\mapsto X+b$ shows that we can suppose that $b=0$, and with $b=0$ the translation $Y\mapsto Y+\sqrt{e}$ gives that we can take $e=0$ as well. Finally $Y\mapsto a^{3}Y$ and $X\mapsto a^{2}X$ gives the form
\[Y^2+XY=X(X^2+AX+B),\qquad A,B\in\Fq\]
as equation for the curve $E$. Hilbert's Satz 90 shows that $\alpha^2+\alpha+A=0$ has a solution $\alpha\in\Fq$ if and only if $\text{Tr}_{\Fq/\F_2}(A)=0$. If this trace equals 1 we can take $\alpha$ in a degree 2 extension of $\Fq$. The change of variables $Y\mapsto Y+\alpha X$ yields then the elliptic curve $E'$ with equation $Y^2+XY=X(X^2+B)$. The conclusion is that $E'\cong E$ over $\Fq$ if $\text{Tr}_{\Fq/\F_2}(A)=0$, and otherwise we have $E'\cong \text{Twist}(E)$. As we did not find a relevant reference we prove the following lemma, which implies that the sum of the traces of Frobenius of $E$ and $E'$ is zero.
\begin{lemma}\label{lem:sumTwistEven}
The equations $Y^2+XY=X(X^2+AX+B)$ and $Y^2+XY=X(X^2+B)$, where $A,B\in\Fq$ and $\text{Tr}_{\Fq/\F_2}(A) = 1$, have together precisely $2q$ affine solutions.\end{lemma}
\textsc{Proof.} We show that for every $x\in\Fq^\times$ one of the equations has two solutions, and the other has none. If $x=0$ both have one solution. Choose $x\in\Fq^\times$. Replacing $Y$ by $xY$ gives as equation $Y^2+Y=x+B/x+A$, and this has (two) solutions if and only if
$\text{Tr}_{\Fq/\F_2}(x+B/x+A)=0$. The linearity of the trace concludes the proof.\ep

Analogously we can find for supersingular curves an equation $Y^2+\og Y=X^3+X^2$ with similar properties as above. We do not work this out, as we do not need it anyway.

%Indeed, if we start with $Y^2+a'Y=X^3+b'X^2+c'X+d'$, we use $Y\mapsto Y+(c'/a') X$, which gives $Y^2+aY=X^3+bX^2+d$. Then we translate by $Y\mapsto Y+\beta$ with $\beta$ a solution of $\beta^2+a\beta+d=0$ --- which could require a twist --- and hence remove $d$ from the equation. So we may suppose that we have $Y^2+aY=X^3+bX^2$. Finally we take $\delta:=\sqrt{b}$ and transform using $Y\mapsto \delta^3Y$ and $X\mapsto \delta^2X$, which together with $\og:=\delta^3a$ gives the result.

Define $H(X) := X$, $Q_f(X,\G) := X^2+\G$ and $\og := B$, then we have proven that it suffices to compute the zeta function of the elliptic curve with equation
\[Y^2+H(X) Y=H(X)Q_f(X,\og),\]
where $H(X),Q_f(X,\G)\in \F_2[X,\G]$ and $\og\in\Fq$. In order to get an elliptic curve for $Y^2+H(X) Y=H(X)Q_f(X,0)$ as well, we only have to translate $\G\mapsto \G+1$, as $Y^2+XY=X(X^2+1)$ does define an elliptic curve. Again, the above transformations can be done very fast in practice. The most time consuming step is computing $\text{Tr}_{\Fq/\F_2}(A)$, which can certainly be done in time $\Ot(n^2)$. The memory requirements are only $\O(n)$. \section{$p$th power Frobenius in odd characteristic}\label{sec:FrobInOddChar}
Now that we have put our elliptic curve --- up to a twist --- in a linear family, we will show how to compute the matrix of the $p$th power Frobenius on its Monsky-Washnitzer cohomology. This cohomology was first considered by Kedlaya in \cite{KedlayaCountingPoints} in an algorithm to count the number of points on hyperelliptic curves in odd characteristic. We have worked out this deformation approach in great detail in \cite{HubrechtsHECOdd}, and we will give a short summary in this section, specified to genus 1 and with $\Fp$ as base field. More details can hence be found in \cite{HubrechtsHECOdd}.

\subsection{A sketch of the deformation theory}\label{ssec:sketchDeformationOdd}
We assume in this section that $p$ is an odd prime. Let $\bar Q(X,\G)\in\Fp[X,\G]$ be of the form explained at the end of section \ref{ssec:OddCharParam}, in particular monic of degree 3 in $X$ and squarefree for $\G=0$. Suppose $\og$ is the parameter such that we need the zeta function of $E:Y^2=\bar Q(X,\og)$, and let the finite field $\Fq=\F_{p^n}$ be defined as $\Fp[x]/\bar\varphi(x)$ with $\bar\varphi(x)$ the minimal polynomial of $\og$ over $\Fp$.
\begin{note}\label{note:extendedDefiningField} The general case can indeed be reduced to this. Suppose that $\og\in\F_{p^m}$ for $1\leq n\leq m$, then \cite{ShoupMinimalPolynomial} shows how to compute the minimal polynomial of $\og$ over $\Fp$ in time $\Ot(m^2)$, and hence also the field $\Fq=\F_{p^n}$. Having computed the zeta function over $\Fq$, we can use formula (\ref{eq:ZetaOverExtension}) to conclude the algorithm.\end{note}
Denote with $\Qp$ the field of $p$-adic numbers, and with $\Qq$ the unique degree $n$ unramified extension of $\Qp$. In fact we need a very specific representation of $\Qq$, which will be explained at the end of section \ref{ssec:computatIssuesOdd}. We write $\Zp$ and $\Zq$ for the rings of integers of these fields, and the Frobenius automorphism, a lift of $x\mapsto x^p$ on $\Fq$, is denoted by $\sigma$. This morphism $\sigma$ is extended with $\sigma(\G):=\G^p$. The valuation on $\Qq$ is denoted by $\ord$, normalized to $\ord(p)=1$.

The Monsky-Washnitzer construction starts with a degree preserving lift $Q(X,\G)\in\Zp[X,\G]$ of $\bar Q(X,\G)$. Define the resultant $$r(\G):=\text{Res}_X\left(Q(X,\G);\frac{\partial}{\partial X}Q(X,\G)\right),$$ then we find that $\bar r(0)$ and $\bar r(\og)$ (where $\bar{\ }$ denotes the reduction modulo $p$) are both nonzero due to the fact that $0$ and $\og$ give (nonsingular) elliptic curves. Write $r(\G)=\sum r_i\G^i$ and let $\rho'$ be the largest index $i$ such that $\ord(r_i)=0$. Then we define $R(\G):=\sum_{i=0}^{\rho'}r_i\G^i$, so that $R(\G)$ has a unit in $\Zp$ as leading coefficient and $R(\G)\equiv r(\G)\bmod p$. Define the ring $S:=\Qp[\G,1/R(\G)]^\dagger$, where $\dagger$ denotes the overconvergent completion, and the $S$-module
\[T := \frac{\Qp[X,Y,1/Y,\G,1/R(\G)]^\dagger}{(Y^2-Q(X,\G))}.\] On $T$ act two differential operators, namely $d:T\to TdX:v\mapsto \frac{\partial v}{\partial X}dX$ and the connection $\nabla:T\to Td\G:v\mapsto \frac{\partial v}{\partial \G}d\G$. The submodule $H_{MW}^-$ of $TdX/dT$ is defined as the eigenspace under the elliptic involution, and is a free 2-dimensional $S$-module. With $F_p$ the Frobenius map on $H_{MW}^-$, we find the following commutative diagram:
\begin{equation}\label{eq:diagram}\begin{CD}
H_{MW}^- @>{\nabla}>> H_{MW}^-d\G\\
@VV{F_p}V @VV{F_p}V\\
H_{MW}^- @>{\nabla}>> H_{MW}^-d\G.
\end{CD}\end{equation}
The basis used in \cite{HubrechtsHECOdd} for $H_{MW}^-$ is the pair $\{dX/\sqrt{Q},XdX/\sqrt{Q}\}$, and the diagram (\ref{eq:diagram}) gives the differential equation
\begin{equation}\label{eq:diffEq}\frac{\partial}{\partial \G} F(\G)+F(\G)G(\G) = p \G^{p-1} G^{\sigma}(\G^{p}) F(\G).\end{equation}
for the matrix $F(\G)$ of $F_p$ with respect to this basis. Here $G(\G)$ is the matrix of the connection $\nabla$. Let $\g$ be the Teichm\"{u}ller lift of $\og$ in $\Zq$, then the matrix $F(\g)$ is precisely the matrix of the $p$th power Frobenius on the Monsky-Washnitzer cohomology as found by Kedlaya in \cite{KedlayaCountingPoints}.\\

\subsection{Computational issues}\label{ssec:computatIssuesOdd}
In section \ref{sec:Eigenvalue} we will need the matrix $F(\g)$ up to a certain $p$-adic precision $N=\O(n)$. Following the algorithm in \cite{HubrechtsHECOdd} with $g=a=\kappa=1$, and limiting ourselves to steps 1 to 7 of the algorithm, this can be achieved in time $\Ot(n^2)$ and space $\O(n^2)$.

There are two important points to note. First, we will need that $F(\g)$ is integral, which is a priori not guaranteed with our chosen basis if $p=3$ (see \cite[section 3.5]{KedlayaComputingZetaFunctions}). Two possible solutions emerge. We can imitate the proofs of \cite{HubrechtsHECOdd}, but now with the basis $\{dX/\sqrt{Q}^3,XdX/\sqrt{Q}^3\}$, which does give an integral matrix. The complexity estimates will all remain the same in this case, and this is the solution used in the implementation we made. Another possible work-around is to compute the matrix of the change between the two bases, a matrix that becomes integral after multiplying with $p$ and is easily retrieved using Kedlaya's algorithm. Transforming $F(\g)$ using this matrix yields then an integral version of $F(\g)$.

Second, in the algorithm a particular representation of $\Qq=\Qp[x]/\varphi(x)$ is used, namely $\varphi(x)$ has to be a \emph{Teichm\"{u}ller modulus} lift of $\bar\varphi(x)$. This means simply that both polynomials are equal modulo $p$, and that $\varphi(x)$ is a monic divisor of $x^q-x$. Equivalently we can say that $\varphi(x)$ is the minimal polynomial of the Teichm\"{u}ller lift $\g$ of $\og$. In \cite[section 12.1.2]{CohenFrey} a very efficient algorithm for computing $\varphi(x)$ is given, originally due to Harley, that computes $\varphi(x)$ in time $\Ot(n^2)$ and space $\O(n^2)$. \section{$2$nd power Frobenius in characteristic 2}\label{sec:FrobInChar2}
As proven in section \ref{ssec:EvenCharParam}, it suffices to consider curves given by\[Y^2+XY=X(X^2+\og+1),\quad\og\in\Fq,\ q= 2^n.\] Again we will explain briefly how to compute the matrix of the second power Frobenius on the Monsky-Washnitzer cohomology of the curve. It was first shown in \cite{DenefVercauteren} how to do this in time $\Ot(n^3)$ and space $\O(n^3)$, and in \cite{HubrechtsHECEven} we extended this result so that it worked faster and used less memory in one dimensional families. We will now sketch how this works, all details can be found in \cite{HubrechtsHECEven}.

\subsection{Computing the matrix of Frobenius}\label{ssec:ComputeFrobEven}
We suppose as in the previous section that $\Fq$ is given by $\F_2[x]$ divided by the minimal polynomial of $\og$. Define $\Q_2$, $\Qq$, $\Z_2$, $\Zq$ and $\sigma$ as before, and let $H(X):= X$ and $Q_f(X,\G) := X^2+\G+1$. The polynomial $c(\G)$ from \cite{HubrechtsHECEven} is just equal to 1 in our case. The resultant needed is $r(\G)=\text{Res}_{X}(H;Q_f\frac{\partial}{\partial X}H)=\G+1$, and clearly both $\bar r(0)$ and $\bar r(\og)$ are nonzero in $\Fq$. Moreover, defining $R(\G)$ as before yields $R(\G)=r(\G)$. The ring $S$ is defined by $S:=\Q_2[\G,1/(\G+1)]^\dagger$ and the $S$-module $T$ by
\[T := \frac{\Q_2[X,Y,1/X,\G,1/(\G+1)]^\dagger}{(Y^2+XY-X(X^2+\G+1))}.\] Using the definitions of $d$, $\nabla$, $H_{MW}^-$ as before, we find again diagram (\ref{eq:diagram}), with $\mathcal{B} := \{YdX,XYdX\}$ as basis for $H_{MW}^-$. Here too we get $F(\g)$, using the Teichm\"{u}ller modulus representation of $\Qq$, with precision $N=\O(n)$. However, in order to get an integral matrix our chosen basis does not suffice, indeed, from the proof of proposition 11 from \cite{HubrechtsHECEven} follows that only $2^6 \cdot F(\g)$ is guaranteed to be integral. We will show in the next subsection how to solve this problem. The conclusion will be that we have to compute $F(\g)$ modulo $2^{N+10}$, and can transform it afterwards into a matrix of Frobenius modulo $2^N$ with integral coefficients. As follows from the algorithm of \cite{HubrechtsHECEven}, we can find this approximation of $F(\g)$ in time $\Ot(n^2)$ and space $\O(n^2)$.

We would like to mention that in \cite{GerkmannEC} Gerkmann considered a deformation for the same family $Y^2+XY=X(X^2+\g)$ that we used above.

\subsection{An integral matrix of Frobenius}\label{ssec:integralFrobEven}
We now give a sketch of how to remedy the `integrality problem'. The eigenvalues of the $q$th power Frobenius map are the reciprocal zeroes of the numerator of the zeta function, and hence integral. This implies that a $\Zq$-submodule of $H_{MW}^-$ does exist that is stable under this map. In \cite[Proposition 5.3.1]{EdixhovenCourse}, Edixhoven showed how to find a basis for this submodule, and in \cite{DVErratum} Denef and Vercauteren applied this to their characteristic 2 situation. It turns out that $\mathcal{D} := \{\frac{dX}{2Y+X},\frac{XdX}{2Y+X}\}$ is such an `integral basis'. It might be possible to reconstruct the algorithm explained above using this basis, but here we will explain how to use the matrix of the change of basis in order to achieve an integral matrix of Frobenius.

We now briefly recall the result of \cite{DVErratum}, specialized to our situation. The modules $H_1$ and $H_1^-$ are as defined in Denef and Vercauteren's paper \cite{DenefVercauteren}, essentially they are the modules $TdX/dT$ and $H_{MW}^-$ above specialized to $\G=\g$. The curve $E:Y^2+XY-X(X^2+\g+1)=0$ is a smooth and proper curve over $\Zq$, and $E\backslash \{P_\infty\}$ is affine, with $P_\infty$ the point at infinity of $E$. Let $D=kP_\infty$ be a divisor on $E$ with $k\geq 2$. We define the $\Zq$-module $L$ as consisting of those differentials $\omega$ on $E\backslash\{P_\infty\}$ satisfying the following two conditions. First, we require that $\text{div}(\omega)+D\geq 0$, and second, each term with valuation less than $-1$ in the local expansion of $\omega$ at $P_\infty$ is integrable over $\Zq$. Then the image of $L$ in $H_1$ is independent of the choice of the divisor $D$ and invariant under the $p$th power Frobenius, and $L$ generates $H_1$. Hence we have also that $L\cap H_1^-$ generates $H_1^-$, and $\mathcal{D}$ will be a basis for $L\cap H_1^-$ as $\Zq$-module.\\

First we need a lower bound on the valuation of the matrix of change of basis and its inverse. The differential forms $YdX$ and $XYdX$ from $\mathcal{B}$ have a pole of order 6 respectively 8 at the point $P_\infty$. If we take $D=8P_\infty$, both forms satisfy the condition $\text{div}(\omega)+D\geq 0$, and $4\omega$ for $\omega\in\mathcal{B}$ will also satisfy the second condition on the integrability. Indeed, during integration only $-7,\ldots,-1$ can appear as denominators, and 4 divided by one of these is always integral in $\Z_2$. This implies that both $4YdX$ and $4XYdX$ are in the $\Zq$-module $L$, which has $\mathcal{D}$ as basis, and hence the matrix defining the change of basis $\mathcal{B}$ to $\mathcal{D}$ has valuation at least $-2$.

For the inverse we have to reduce the basis $\mathcal{D}$ to $\mathcal{B}$ and use the lemmata 2 and 3 of \cite{DenefVercauteren}. As
\[\frac{dX}{2Y+X} = \frac{(2Y+X)dX}{4X(X^2+\g+1)+X^2}= \frac{(2Y+X)}{X^2}dX\cdot\left(\sum_{k=0}^\infty (-4)^k\left(X+\frac{\g+1}{X}\right)^k\right),\]
an easy computation gives as lower bound for the matrix of this change of basis
\[\min\left\{ \min_{k\geq 2}(2k-3-\lfloor\log_2(k+1)\rfloor); \min_{k\geq 1}(2k-3-\lfloor\log_2(k+3)\rfloor) \right\}+1\geq -2.\]
Computing this last matrix, denoted with $B$, modulo $2^M$ with $M=\O(n)$ is easy using the reduction formulae in \cite{DenefVercauteren}, but this would require time $\Ot(n^3)$. We can see however that we do not need $B$ modulo such high power of $2$. Indeed, let $B'$ be any invertible matrix over $\Qq$ such that $F':=(B'^{-1})^\sigma F(\g)B'$ is integral, then $B'$ gives the change to an unknown but irrelevant basis, and the resulting integral matrix $F'$ is still a matrix of Frobenius. So let $B'\equiv B\bmod 2^\alpha$ for some $\alpha$, then if $(B'^{-1})^\sigma F(\g) B'$ is integral we are done. We will show that $\alpha = \O(1)$ suffices, and as a consequence the algorithm of \cite{DenefVercauteren} allows us to compute $B'$ in time $\Ot(n^2)$ and space $\O(n^2)$.

From the valuation bound $-2$ on $B^{-1}$ above we see that $\ord(\det B)$ is bounded by 4, and hence working modulo $2^6$ suffices already to be able to invert $B'$ (which has to be done to the full precision $2^{N+10}$). As soon as the $2$-adic precision of $B'$ exceeds 8 bits, the sum of the valuation bounds for $F(\g)$ and $B'^{-1}$, the resulting product $(B'^{-1})^\sigma F(\g) B'$ will be an integral matrix of Frobenius as required. Hence taking $\alpha := 9=\O(1)$ suffices. The loss in precision in this product is at most $2+6+2$, hence we have to compute $F(\g)$ modulo $2^{N+10}$. \section{An eigenvalue of the $q$th power Frobenius}\label{sec:Eigenvalue}
In this section we will first show that it suffices to compute an approximation of an eigenvalue of the matrix of the $q$th power Frobenius, and we reduce this to computing an `eigenvalue' of $F(\g)$, in fact an eigenvalue of the $\sigma$-linear Frobenius map $F_p$. In a second subsection we explain how to solve this last problem, by showing that we can always satisfy certain conditions required for an algorithm that computes solutions of a specific type of $p$-adic equation.

\subsection{Reduction to an `eigenvalue' of $F(\g)$}\label{ssec:eigenvalue1}
Suppose that $E$ is a nonsupersingular curve over $\Fq$, where $q= p^n$, and $F=F(\g)$ is the matrix of the $p$th power Frobenius on its Monsky-Washnitzer cohomology over $\Zq$, as explained in the two previous sections. For $$\mathcal{F}:=F^{\sigma^{n-1}}\cdot F^{\sigma^{n-2}}\cdots F^\sigma\cdot F,$$ the matrix of the $q$th power Frobenius, Kedlaya \cite{KedlayaCountingPoints} and Denef and Vercauteren \cite{DenefVercauteren} showed that we have, with $Z(T)$ the zeta function of $E$ over $\Fq$,
\[Z(T) = \frac{\det(1-\mathcal{F}T)}{(1-T)(1-qT)}.\] If we write $qT^2-tT+1$ for the numerator of the zeta function, it follows immediately that $\det(\mathcal{F})=q$ and $\text{Tr}(\mathcal{F})=t$. Let $\lambda_1$ and $\lambda_2$ be the eigenvalues of $\mathcal{F}$, then we will prove in the next subsection that $\lambda_1,\lambda_2\in\Zq$, and that we may suppose that $\ord(\lambda_1)=0$ and hence $\ord(\lambda_2)=\ord(q/\lambda_1)=n$. We are trying to compute $t=\text{Tr}(\mathcal{F})=\lambda_1+q/\lambda_1$. The Hasse-Weil bound says\footnote{In fact, the Hasse-Weil bound shows that $|t|\leq 2\sqrt{q}$, but equality can only occur for supersingular curves.} that $|t|<2\sqrt{q}$, hence we only need to compute $\lambda_1$ modulo $p^N$ with
\begin{equation}\label{eq:DefinitionN}N := \lceil\log_p(4\sqrt q)\rceil= \lceil n/2 + \log_p(4)\rceil = \O(n),\end{equation} which is smaller than $n$ if $n$ is not too small. To conclude, if suffices to compute $\lambda_1$ modulo $p^N$ in order to find the zeta function of $E$: the trace $t$ is then the unique integer congruent to $\lambda$ modulo $p^N$ that satisfies $|t|< 2\sqrt q$.\\

If we have matrices $C$ and $D$ over $\Zq$ such that $F=C^\sigma D C^{-1}$ with $D$ in uppertriangular form, this implies
\[\mathcal{F} = C\cdot \left(D^{\sigma^{n-1}}\cdot D^{\sigma^{n-2}}\cdots D^{\sigma}\cdot D\right)\cdot C^{-1},\]
and with $\mu$ the upper diagonal element of $D$ this gives that the norm $\mathcal{N}_{\Qq/\Qp}(\mu)$ is an eigenvalue of $\mathcal{F}$. We will show in section \ref{ssec:eigenvalue2} that such $\mu$ with valuation 0 can always be found efficiently if $E$ is not supersingular. It is easily seen that a factorization $F=C^\sigma D C^{-1}$ over $\Qq$ cannot exist if the curve is supersingular: the product of the two diagonal elements has valuation one, and their sum has then valuation at least one. This is clearly impossible as the valuation is a map from $\Qq$ to the integers.

Having found $\mu$ we still have to compute its norm. For this we can apply an algorithm by Harley, which uses an adaptation of Moenck's extended \textsc{gcd} algorithm in order to compute a certain resultant. Indeed, if $\Zq=\Zp[x]/\varphi(x)$ with $\varphi(x)$ a Teichm\"{u}ller modulus, and $\mu(x)\in\Zp[x]/\varphi(x)$, then $$\mathcal{N}_{\Qq/\Qp}(\mu) = \text{Res}_x(\mu(x),\varphi(x)).$$ A complete description of the algorithm has been given by Vercauteren and can be found in \cite[Section 3.10.3]{VercauterenThesis}. It requires $\Ot(n^2)$ time and $\O(n^2)$ space. As noted there, in order for the algorithm to work well $\mu(x)$ should have a unit in $\Zp$ as leading coefficient. This is however easily forced: suppose $\mu(x)\bmod p$ has degree $n-1-r$, then $x^r\mu(x)$ satisfies this condition. Moreover, $x^r$ itself satisfies the condition as well, hence computing $\mathcal{N}(\mu)=\mathcal{N}(x^r\mu(x))/\mathcal{N}(x^r)$ gives the required result. Note that $x^r$ is a Teichm\"{u}ller lift with $\mathcal{N}(x^r)=\varphi(0)^r$, and its norm can thus be computed much faster.

\subsection{Computation of an `eigenvalue' $\mu$ of $F(\g)$}\label{ssec:eigenvalue2}
In this subsection $\equiv$ will always mean `congruence modulo $p$', unless `$\bmod\,p^N$' is explicitly written. We will need the following algorithm of Harley, which can be found as algorithm 3.10.2 in \cite{CohenFrey}. Note that this algorithm requires $\Zq$ to be given as $\Zp[x]$ modulo a Teichm\"{u}ller modulus.

\noindent\texttt{INPUT:\hphantom{ }} $\psi(X,Y)\in\Zq[X,Y]$, $x_0\in\Zq$ such that
\[\psi(x_0,x_0^\sigma)\equiv \frac{\partial\psi}{\partial X}(x_0,x_0^\sigma)\equiv 0,\qquad \frac{\partial\psi}{\partial Y}(x_0,x_0^\sigma)\not \equiv 0,\]
\noindent\texttt{OUTPUT:} $\alpha\in\Zq$ such that
\[\alpha\equiv x_0,\qquad \psi(\alpha,\alpha^\sigma)\equiv 0\bmod p^N.\]
Following the complexity estimates found in \cite{CohenFrey}, it is easily shown that if the degree of $\psi$ is fixed, the algorithm runs in time $\Ot(nN)$ and space $\O(nN)$.\\

Write $F=F(\g)$ as $\begin{pmatrix}f_1 & f_2\\f_3 & f_4\end{pmatrix}$ with all $f_i$ in $\Zq$, and consider the system of equations
\begin{equation}\label{eq:Eigenvector}
\begin{pmatrix}f_1 & f_2\\f_3 & f_4\end{pmatrix}
\begin{pmatrix}1 \\ \alpha         \end{pmatrix} = \mu
\begin{pmatrix}1 \\ \alpha^\sigma  \end{pmatrix},\quad \text{or}\quad
\begin{cases}\ f_1+\alpha f_2=\mu,\\\
f_3+\alpha f_4=\mu\alpha^\sigma.\end{cases}\end{equation}
It is clear that if we can find a solution $(\alpha,\mu)\in\Zq\times\Zq^\times \bmod p^N$ for (\ref{eq:Eigenvector}), this yields a factorization of $F=(C^\sigma) D C^{-1}$, which is of the kind that we are looking for. Here $C$ and $D$ equal \[C = \begin{pmatrix}1 & *\\\alpha & *\end{pmatrix},\qquad  D=\begin{pmatrix}\mu & *\\0 & *\end{pmatrix},\qquad \text{all $*\in \Zp$}.\]

Eliminating $\mu$ from the equations (\ref{eq:Eigenvector}) gives
\begin{equation}\label{eq:eqEigenvector}\alpha(\alpha^\sigma f_2-f_4)+ (\alpha^\sigma f_1-f_3)=0.\end{equation}
If $f_1\equiv f_2\equiv 0$, certainly one of $f_3$, $f_4$ will not be zero modulo $p$, as $\ord(\det(F)) = 1$. In this case we can work with $(\alpha\ \ 1)^T$ instead of $(1\ \  \alpha)^T$. So we can suppose that at least one of $f_1$ or $f_2$ is nonzero modulo $p$. Let
\[x_0^\sigma := \left(\frac{f_4}{f_2}\bmod p\right),\qquad \text{or}\qquad x_0^\sigma := \left(\frac{f_3}{f_1}\bmod p\right).\] If both definitions make sense, $\det(F)\equiv 0$ implies that they are equal. Computing the corresponding $x_0$ is easy finite field arithmetic. We define the polynomial $\psi(X,Y)$ by
\[\psi(X,Y) := X(Yf_2-f_4)+(Yf_1-f_3)\ \ \in\Zq[X,Y].\]
Our choice of $x_0^\sigma$ guarantees that $\psi(x_0,x_0^\sigma)\equiv 0$ and also
\[\frac{\partial}{\partial X}\psi(x_0,x_0^\sigma)=x_0^\sigma f_2-f_4\equiv 0.\]
This last inequality holds even if $f_2\equiv 0$.
%Indeed, $f_1$ is then nonzero.
We will show immediately that $\frac{\partial}{\partial Y}\psi(x_0,x_0^\sigma)\not\equiv 0$ follows from nonsupersingularity. The algorithm from the beginning of this section allows us now to compute $\alpha\in\Zq$ and hence $\mu$ with precision $N=\O(n)$ in time $\Ot(n^2)$ and $\O(n^2)$. In addition, eliminating $\alpha$ from (\ref{eq:Eigenvector}) yields
\[\mu(f_4-\alpha^\sigma f_2)=f_1 f_4 - f_2f_3,\]
which equals $\det(F)$ and has valuation 1. As $f_4-\alpha^\sigma f_2\equiv 0$, it is impossible that $\ord(\mu)>0$ as well.\\

Suppose that $\frac{\partial}{\partial Y}\psi(x_0,x_0^\sigma)=x_0f_2+f_1\equiv 0$. If $f_2\equiv 0$ this would imply $f_1\equiv 0$, which we excluded. Define $f_i'=f_i/f_2$, then $$f_1'\equiv -x_0,\quad f_4'\equiv x_0^\sigma\equiv x_0^p,\quad f_3'\equiv f_1'f_4'\equiv -x_0^{p+1}.$$ As a consequence $$F\equiv f_2\begin{pmatrix}-x_0&1\\-x_0^{p+1} & x_0^p \end{pmatrix}\qquad\text{and}\qquad F^\sigma F\equiv \begin{pmatrix}0&0\\0&0\end{pmatrix}.$$ This implies that the trace of $\mathcal{F}$ is congruent to zero modulo $p$, and hence the curve considered is supersingular. \section{Conclusion and implementation results}\label{sec:conclusion}
Combining all steps explained in sections \ref{sec:parameterFamily}, \ref{sec:FrobInOddChar}, \ref{sec:FrobInChar2} and \ref{sec:Eigenvalue} above, we have found a deterministic algorithm that for every elliptic curve over $\F_{p^n}$ given by its Weierstrass equation, can compute its zeta function in time $\Ot(n^2)$ and space $\O(n^2)$. We will now give a list of the main steps of the algorithm. We assume that we are working in odd characteristic, and with an `integral basis' for the Monsky-Washnitzer cohomology $H_{MW}^-$. We do not mention in the algorithm that we only compute \emph{approximations} of the objects involved. If $p^n$ is so small that $N>n$ in (\ref{eq:DefinitionN}), we can use a naive point counting algorithm.\\

\noindent\texttt{INPUT:\hphantom{\ }} Finite field $\F_{p^n}$, monic squarefree polynomial $Q(X)\in\F_{p^n}[X]$ of degree 3,\vspace{.5 \baselineskip}

\noindent\texttt{OUTPUT:} The zeta function of the elliptic curve $Y^2=Q(X)$ over $\F_{p^n}$.\vspace{.5 \baselineskip}

\noindent\texttt{STEP 1:} Put the curve in a one parameter family $Y^2=Q(X,\og)$, where $\og\in\F_{p^n}$, as explained in section \ref{sec:parameterFamily}. \vspace{.5 \baselineskip}

\noindent\texttt{STEP 2:} Compute the matrix of Frobenius $F(0)$ of $Y^2=Q(X,0)$, and the differential equation for $F(\G)$.\vspace{.5 \baselineskip}

\noindent\texttt{STEP 3:} Solve the differential equation and find $F(\G)\in\Zp[[\G]]^{2\times 2}$.\vspace{.5 \baselineskip}

\noindent\texttt{STEP 4:} Determine $\bar\varphi(x)\in\Fp[x]$, the minimal polynomial of $\og$, and define $\F_{p^m} := \Fp[x]/\bar\varphi(x)$. Lift $\bar \varphi(x)$ to a Teichm\"{u}ller modulus so that $\Z_{p^m} = \Zp[x]/\varphi(x)$ and $x=\g$.\vspace{.5 \baselineskip}

\noindent\texttt{STEP 5:} Compute $F(\g)$ by reducing $F(\G)$ modulo $\varphi(\G)$.\vspace{.5 \baselineskip}

\noindent\texttt{STEP 6:} Compute a solution $(\alpha,\mu)$ with $\ord(\mu)=0$ for the equation
\[F(\g)\cdot
\begin{pmatrix}1 \\ \alpha         \end{pmatrix} = \mu
\begin{pmatrix}1 \\ \alpha^\sigma  \end{pmatrix}.\]

\noindent\texttt{STEP 7:} Compute $t_1\equiv \mathcal{N}_{\Q_{p^m}/\Qp}(\mu)$ modulo an appropriate power of $p$, such that $|t_1|<2\sqrt{p^m}$. Compute then the resultant $$p^nT^2-tT+1 = \text{Res}_X(p^mX^2-t_1X+1,X^{n/m}-T).$$

\noindent\texttt{STEP 8:} Return
\[\frac{p^nT^2-tT+1}{(1-T)(1-p^nT)}.\]\\

\noindent We have implemented this algorithm in odd characteristic, and present a few timing results obtained with it. Note that we do \emph{not} use Harley's $\Ot(n^2)$ norm algorithm for step 7, but instead the --- far easier to implement and in practice probably faster for reasonable $n$ --- algorithm of Satoh, Skjernaa and Taguchi \cite{sst}. This method runs in time $\Ot(n^{2.5})$ given some precomputations. These precomputations require time $\Ot(n^3)$, but are completely integer arithmetic and hence extremely fast. In our algorithm they are necessarily part of the algorithm (they depend on $\bar\varphi(x)$, the minimal polynomial of the parameter $\og$), so our implementation has as theoretical complexity $\Ot(n^3)$. In step 2 the matrix $F(0)$ is computed using an implementation of Kedlaya's algorithm by Michael Harrison.

The implementation has been made in the computational algebra system Magma V2.13-3, and the timing results were obtained on an AMD Athlon 64 3000+, using 1GB of physical memory. The algorithm received as input a random elliptic curve over $\F_{p^n}$, given by its Weierstrass equation. All times in the following table are in seconds.

\begin{center}{\small
\begin{tabular}{|r||r|r|r|r|r|r|r|}
\hline $p\backslash n$\vphantom{$\sum^j$} & 50 & 100 & 250 & 500 & 1000 & 2000 & 4000\\\hline\hline
3\vphantom{$\sum^j$} & .35  & .78  & 4.47  & 18.43 & 99 & 604 & 4293\\\hline
5\vphantom{$\sum^j$} & 1.18 & 2.73 & 11.43 & 47.92 & 227 & 1389 & -\\\hline
7\vphantom{$\sum^j$} & 3.27 & 7.79 & 41.85 & 186.48 & 957 & 5592 & -\\\hline
\end{tabular}}
\end{center}
It is interesting to see that for $n\gg 0$ almost all computation time goes to steps 6, 7 and the computation of the Teichm\"{u}ller modulus in step 4, the first two being comparable in required time. E.g.\ for $p^n=3^{4000}$ we have as total time 4293 seconds, where step 6 uses 1916 seconds and step 7 uses 2190 seconds. For $p^n=7^{2000}$ the computation of $\varphi(x)$ takes 3910 seconds. A conclusion that could be drawn from this is that for such big fields our algorithm should work faster than Harley's --- as long as in either algorithm the same norm algorithm and no precomputation is used --- because he needs a computation similar to step 6 but with an equation $\psi$ of higher degree, and exactly the same field polynomial and norm computation.

Steps 2 and 3 can be considered as precomputation, meaning that they only depend on the field size (and the structure of the family in which the curve lives). For $n$ big enough these steps are of minor influence, but for fields of cryptographic size it is worth looking at the time needed for just one curve. The following table gives these times for the field sizes as above, hence ignoring the time for steps 2 and 3 of the algorithm.
\begin{center}{\small
\begin{tabular}{|r||r|r|r|r|r|r|r|}
\hline $p\backslash n$\vphantom{$\sum^j$} & 50 & 100 & 250 & 500 & 1000 & 2000 & 4000\\\hline\hline
3\vphantom{$\sum^j$} & .15  & .46  & 3.65 & 16.58 & 95 & 592 & 4252\\\hline
5\vphantom{$\sum^j$} & .26  & .92 & 6.77 & 36.55 & 198 & 1306 & -\\\hline
7\vphantom{$\sum^j$} & 1.76 & 4.87 & 34.06 & 167.56 & 909 & 5447 & -\\\hline
\end{tabular}}
\end{center} 
\fontsize{9}{11pt}\selectfont
\bibliographystyle{acm}
\bibliography{bibliography}

\end{document}